\documentstyle[12pt,amstex,amscd,psfig]{amsart}

\begin{document}

\newtheorem{theorem}{Theorem}[section]
\newtheorem{prop}[theorem]{Proposition}
\newtheorem{lemma}[theorem]{Lemma}
\newtheorem{cor}[theorem]{Corollary}

\title{Cannon-Thurston Maps for Trees of Hyperbolic Metric Spaces}

\author{Mahan Mitra}\thanks{Research partly supported by Alfred P. Sloan
Doctoral Dissertation Fellowship, Grant No. DD 595. }

\begin{abstract} 
Let (X,d) be a tree (T) of hyperbolic metric spaces satisfying the
quasi-isometrically embedded condition.  Let $v$ be a vertex of $T$. Let
$({X_v},d_v)$ denote the hyperbolic metric space corresponding to $v$.
 Then $i : X_v \rightarrow X$ extends continuously to a map
$\hat{i} : \widehat{X_v} \rightarrow \widehat{X}$. This generalizes a Theorem
of Cannon and Thurston. The techniques are
 used to give a new proof of a result of
Minsky: Thurston's 
ending lamination conjecture for certain Kleinian groups. Applications to
graphs of hyperbolic groups and local connectivity of limit sets of Kleinian
groups are also given.
\end{abstract}

\maketitle

\section{Introduction}

Let $G$ be a hyperbolic group in the sense of Gromov \cite{Gromov}. 
Let $H$ be a hyperbolic subgroup of  $G$. We
choose a finite symmetric generating set for $H$ and extend it to a finite
symmetric generating set for $G$. 
Let $\Gamma_H$ and $\Gamma_G$ denote the
Cayley graphs of $H$, $G$ respectively with respect to these generating sets.
By 
adjoining the Gromov boundaries $\partial\Gamma_H$ and $\partial\Gamma_G$
 to $\Gamma_H$ and $\Gamma_G$, one obtains their compactifications
$\widehat{\Gamma_H}$ and $\widehat{\Gamma_G}$ respectively.

We'd like to understand the extrinsic geometry of $H$ in $G$. Since the
objects of study here come under the purview of coarse geometry, asymptotic
or `large-scale' information is of crucial importance. That is to say, one
would like to know what happens `at infinity'.  We put this in the more
general context of a hyperbolic group $H$ acting freely and properly discontinuously
on a proper hyperbolic metric space $X$. Then there is a natural map
 $i : \Gamma_H \rightarrow X$, sending  the vertex set of $\Gamma_H$ to
the orbit of a point under $H$, and connecting images of adjacent vertices
in $\Gamma_H$ by geodesics in $X$. Let $\widehat{X}$ denote the Gromov
compactification of $X$.
 
 A natural question seems to be the following:

\smallskip

{\bf Question:} Does 
the continuous proper map $ i$ : $\Gamma_H \rightarrow X$
extend to a continuous map $\hat i$ : 
$\widehat{\Gamma_H} \rightarrow \widehat{X}$ ?

\medskip

Questions along this line have been raised by Bonahon \cite{Bonahon}.
Related questions in the context of Kleinian groups have been studied by
Bonahon \cite{Bonahon1}, Floyd \cite{Floyd} and  Minsky \cite{Minsky}.
In, \cite{CT}, \cite{Floyd} or \cite{Minsky}, explicit metrics were used.
So though some of their results can be thought of as `coarse', the 
techniques of proof are
not. In \cite{mitra1}, coarse techniques were used to answer the above
question affirmatively for $X = \Gamma_G$, where $G$ is a hyperbolic
group and $H$ a normal subgroup of $G$. This in turn was a generalization of
 a theorem of Cannon and Thurston \cite{CT}. In this paper,
we extend the techniques of \cite{mitra1} to cover examples
arising from  trees of hyperbolic metric spaces  satisfying
an extra technical condition introduced by Bestvina and Feighn in 
\cite{BF}:  the {\it quasi-isometrically embedded} condition. [See Section 3 of
this paper or \cite{BF} for definitions.]

{\bf Definition:} { \it Let $X$ and $Y$ be hyperbolic metric spaces and
$i : Y \rightarrow X$ be an embedding. 
 A {\bf Cannon-Thurston map} $\hat{i}$  from $\widehat{Y}$ to
 $\widehat{X}$ is a continuous extension of $i$. Such a continuous extension
will occassionally be called a Cannon-Thurston map for the pair $(Y,X)$. 
 If $Y = \Gamma_H$ and $X = \Gamma_G$ for a
hyperbolic subgroup $H$ of a hyperbolic group $G$, a Cannon-Thurston map
for $({\Gamma_H},{\Gamma_G})$ will occassionally be referred to as a Cannon-Thurston
map for $(H,G)$.
}

It is easy to see that such a continuous extension, if it exists, is unique.

\medskip
The main theorem of this paper is :

\medskip

{\bf Theorem \ref{main}:}
{\it Let (X,d) be a tree (T) of hyperbolic metric spaces satisfying the
quasi-isometrically embedded condition.  Let $v$ be a vertex of $T$. Let
$({X_v},d_v)$ denote the hyperbolic metric space corresponding to $v$.
 Then $i : X_v \rightarrow X$ extends continuously to a map
$\hat{i} : \widehat{X_v} \rightarrow \widehat{X}$.}

\medskip

A direct consequence of Theorem \ref{main} above is the following:

{\bf Corollary \ref{main1}:}
{\it Let $G$ be a hyperbolic group acting cocompactly on a simplicial tree
$T$ such that all vertex and edge stabilizers are hyperbolic. Also
suppose that every inclusion of an edge stabilizer in a vertex stabilizer
is a quasi-isometric embedding. Let $H$ be the stabilizer of a vertex or
edge of $T$. Then 
there exists a Cannon-Thurston map from
$\widehat{\Gamma_H}$ to $\widehat{\Gamma_G}$.}

\medskip

In \cite{BF}, Bestvina and Feighn give sufficient conditions for
a graph of hyperbolic groups to be hyperbolic. Vertex and edge
subgroups are thus natural examples of hyperbolic subgroups of hyperbolic
groups. 
Essentially all previously 
known examples of non-quasiconvex hyperbolic subgroups 
of hyperbolic
groups arise this way. Theorem \ref{main} shows that these have
Cannon-Thurston maps.

Another consequence of Theorem \ref{main} above is:

{\bf Theorem \ref{main2}:} Let $\Gamma$ be a geometrically tame Kleinian 
group, such that ${{\Bbb{H}}^3}/{\Gamma} = M$ has injectivity radius 
uniformly bounded below by some $\epsilon > 0$. Then there exists a continuous
map from the Gromov boundary of $\Gamma$ (regarded as an abstract group)
to the limit set of $\Gamma$ in ${\Bbb{S}}^2_{\infty}$.

Theorem \ref{main2} 
was independently proven by Klarreich \cite{klarreich}, using
different techniques. 

After some further work and using a theorem of Minsky \cite{Minsky2},
we are able to give a different proof of another result of 
Minsky \cite{Minsky} :  Thurston's Ending Lamination Conjecture
for geometrically tame manifolds with 
freely indecomposable fundamental group and a uniform lower bound on
injectivity radius.

{\bf Theorem \ref{main3} }\cite{Minsky}: {\it 
Let $N_1$ and $N_2$ be homeomorphic
hyperbolic 3-manifolds with freely indecomposable fundamental group. Suppose
there exists a uniform
lower bound $\epsilon > 0$  on the injectivity radii of $N_1$ and $N_2$.
If the end invariants of corresponding ends of $N_1$ and $N_2$ are
equal, then $N_1$ and $N_2$ are isometric.}

\section{Preliminaries}

We start off with some preliminaries about hyperbolic metric
spaces  in the sense
of Gromov \cite{Gromov}. For details, see \cite{CDP}, \cite{GhH}. Let $(X,d)$
be a hyperbolic metric space. The 
{\bf Gromov boundary}
 $X$, denoted by $\partial{X}$,
is the collection of equivalence classes of geodesic rays $r:[0,\infty)
\rightarrow\Gamma$ with $r(0)=x_0$ for some fixed ${x_0}\in{X}$,
where rays $r_1$
and $r_2$ are equivalent if $sup\{ d(r_1(t),r_2(t))\}<\infty$.
Let $\widehat{X}$=$X\cup\partial{X}$ denote the natural
 compactification of $X$ topologized the usual way(cf.\cite{GhH} pg. 124).

The {\bf Gromov inner product}
 of elements $a$ and $b$ relative to $c$ is defined 
by 
\begin{center}
$(a,b)_c$=1/2$[d(a,c)+d(b,c)-d(a,b)]$.
\end{center}

\medskip

{\bf Definitions:} {\it A subset $X$ of $\Gamma$ is said to be 
{\bf $k$-quasiconvex}
if any geodesic joining $a,b\in X$ lies in a $k$-neighborhood of $X$.
A subset $X$ is {\bf quasiconvex} if it is $k$-quasiconvex for some $k$.
A map $f$ from one metric space $(Y,{d_Y})$ into another metric space 
$(Z,{d_Z})$ is said to be
 a {\bf $(K,\epsilon)$-quasi-isometric embedding} if
 
\begin{center}
${\frac{1}{K}}({d_Y}({y_1},{y_2}))-\epsilon\leq{d_Z}(f({y_1}),f({y_2}))\leq{K}{d_Y}({y_1},{y_2})+\epsilon$
\end{center}
If  $f$ is a quasi-isometric embedding, 
 and every point of $Z$ lies at a uniformly bounded distance
from some $f(y)$ then $f$ is said to be a {\bf quasi-isometry}.
A $(K,{\epsilon})$-quasi-isometric embedding that is a quasi-isometry
will be called a $(K,{\epsilon})$-quasi-isometry.

A {\bf $(K,\epsilon)$-quasigeodesic}
 is a $(K,\epsilon)$-quasi-isometric embedding
of
a closed interval in $\Bbb{R}$. A $(K,0)$-quasigeodesic will also be called
a $K$-quasigeodesic.}

\medskip

Let $(X,{d_X})$ be a hyperbolic metric space and $Y$ be a subspace that is
hyperbolic with the inherited path metric $d_Y$.
By 
adjoining the Gromov boundaries $\partial{X}$ and $\partial{Y}$
 to $X$ and $Y$, one obtains their compactifications
$\widehat{X}$ and $\widehat{Y}$ respectively.

Let $ i :Y \rightarrow X$ denote inclusion.

\medskip

{\bf Definition:} { \it Let $X$ and $Y$ be hyperbolic metric spaces and
$i : Y \rightarrow X$ be an embedding. 
 A {\bf Cannon-Thurston map} $\hat{i}$  from $\widehat{Y}$ to
 $\widehat{X}$ is a continuous extension of $i$. }

\medskip

The following  lemma
 says that a Cannon-Thurston map exists
if for all $M > 0$ and $y{\in}Y$, there exists $N > 0$ such that if $\lambda$
lies outside an $N$ ball around $y$ in $Y$ then
any geodesic in $X$ joining the end-points of $\lambda$ lies
outside the $M$ ball around $i(y)$ in $X$.
For convenience of use later on, we state this somewhat
differently. The proof is similar to Lemma 2.1 of \cite{mitra1}

\begin{lemma}
A Cannon-Thurston map from $\widehat{Y}$ to $\widehat{X}$
 exists if  the following condition is satisfied:

There exists a non-negative function  $M(N)$, such that 
 $M(N)\rightarrow\infty$ as $N\rightarrow\infty$ and for all geodesic segments
 $\lambda$  lying outside an $N$-ball
around ${y_0}\in{Y}$  any geodesic segment in $\Gamma_G$ joining
the end-points of $i(\lambda)$ lies outside the $M(N)$-ball around 
$i({y_0})\in{X}$.

\label{contlemma}
\end{lemma}
{\it Proof:}
 Suppose $i:Y\rightarrow{X}$
does not extend continuously . Since $i$ is proper, there exist 
sequences $x_m$, $y_m$ $\in{Y}$ and $p\in\partial{Y}$,
such that $x_m\rightarrow p$
and $y_m\rightarrow p$ in $\widehat{Y}$, but $i(x_m)\rightarrow u$
and $i(y_m)\rightarrow v$ in $\widehat{X}$, where 
$u,v\in\partial{X}$ and $u\neq v$.

Since $x_m\rightarrow p$ and $y_m\rightarrow p$, any geodesic in $Y$
joining $x_m$ and $y_m$ lies outside an $N_m$-ball ${y_0}\in{Y}$,
 where $N_m\rightarrow\infty$ as $m\rightarrow\infty$. Any
bi-infinite geodesic in $X$ joining  $u,v\in\partial{X}$
has to pass through some $M$-ball around $i({y_0})$ in $X$ as
$u\neq v$. There exist constants $c$ and $L$ such that for all $m > L$
any geodesic joining $i(x_m)$ and $i(y_m)$ in $X$ 
passes through an $(M+c)$-neighborhood
 of $i({y_0})$. 
Since $(M+c)$ is a constant not depending on the index $m$
this proves the lemma. $\Box$

\medskip

\section{Trees of Hyperbolic Metric Spaces}

We start with a notion closely related to one  introduced in \cite{BF}.

{\bf Definition:} {\it A  tree (T) of hyperbolic metric spaces satisfying
the q(uasi) i(sometrically) embedded condition is a metric space $(X,d)$
admitting a map $P : X \rightarrow T$ onto a simplicial tree $T$, such
that there exist $\delta{,} \epsilon$ and $K > 0$ satisfying the following: \\
\begin{enumerate}
\item  For all vertices $v\in{T}$, 
$X_v = P^{-1}(v) \subset X$ with the induced path metric $d_v$ is a 
$\delta$-hyperbolic metric space. Further, the
inclusions ${i_v}:{X_v}\rightarrow{X}$ 
are uniformly proper, i.e. for all $M > 0$, $v\in{T}$ and $x, y\in{X_v}$,
there exists $N > 0$ such that $d({i_v}(x),{i_v}(y)) \leq M$ implies
${d_v}(x,y) \leq N$.
\item Let $e$ be an edge of $T$ with initial and final vertices $v_1$ and
$v_2$ respectively.
Let $X_e$ be the pre-image under $P$ of the mid-point of  $e$.  
Then $X_e$ with the induced path metric is $\delta$-hyperbolic.
\item There exist maps ${f_e}:{X_e}{\times}[0,1]\rightarrow{X}$, such that
$f_e{|}_{{X_e}{\times}(0,1)}$ is an isometry onto the pre-image of the
interior of $e$ equipped with the path metric.
\item ${f_e}|_{{X_e}{\times}\{{0}\}}$ and 
${f_e}|_{{X_e}{\times}\{{1}\}}$ are $(K,{\epsilon})$-quasi-isometric
embeddings into $X_{v_1}$ and $X_{v_2}$ respectively.
${f_e}|_{{X_e}{\times}\{{0}\}}$ and 
${f_e}|_{{X_e}{\times}\{{1}\}}$ will occassionally be referred to as
$f_{v_1}$ and $f_{v_2}$ respectively.
\end{enumerate}   }

$d_v$ and $d_e$ will denote path metrics on $X_v$ and $X_e$ respectively.
$i_v$, $i_e$ will denote inclusion of $X_v$, $X_e$ respectively into $X$.

The main theorem of this paper can now be stated:

{\bf Theorem: \ref{main}}
 Let (X,d) be a tree (T) of hyperbolic metric spaces satisfying the
qi-embedded condition.  Let $v$ be a vertex of $T$. Then 
$i_v : X_v \rightarrow X$ extends continuously to ${\hat{i_v}} : \widehat{X_v}
\rightarrow \widehat{X}$.

Some aspects of
the proof of the main theorem of this section are similar to the proof of the
main theorem of \cite{mitra1}. 
Given  a geodesic segment $\lambda\subset{X_v}$, we construct a quasi-convex
set
$B_{\lambda}\subset{X}$.
It follows from the construction that
  if $\lambda$ lies outside a large ball around 
$y_0 \in X_v$, 
$B_{\lambda}$ lies outside a large ball around ${i_v}({y_0}) \in X$,
i.e.  for all $M\geq{0}$ there exists $N\geq{0}$
such that if $\lambda$ lies oustside the $N$-ball around $y_0 \in X_v$, 
$B_{\lambda}$ lies outside the $M$-ball around ${i_v}({y_0}) \in X$.
Combining this with Lemma \ref{contlemma} above, the proof of Theorem
\ref{main} is completed.

\medskip

For convenience of exposition, $T$ shall be assumed to be rooted, i.e.
equipped with a base vertex $v_0$. Let $v \neq v_0$ be a vertex of $T$.
Let $v_{-}$ be the penultimate vertex on the geodesic edge path
from $v_0$ to $v$. Let $e$ denote the directed edge from ${v_{-}}$
to $v$.  Define 
${\phi_v} : {f_{e_{-}}}({X_{e_{-}}}{\times}\{{0}\}) \rightarrow
{f_{e_{-}}}({X_{e_{-}}}{\times}\{{1}\})$ by
\begin{center}
${\phi_v}({f_{e_{-}}}({x}{\times}\{{0}\}) =
{f_{e_{-}}}({x}{\times}\{{1}\})$  for $x \in {X_{e_{-}}}$.
\end{center}

Let $\mu$  be a geodesic in 
$(X_{v_{-}})$, joining 
$ {i_{v_{-}}^{-1}}(a), {i_{v_{-}}^{-1}}(b)
 \in 
{i_{v_{-}}^{-1}}{\cdot}{f_{e_{-}}}({X_{e_{-}}}{\times}\{{0}\}) $. ${\Phi}_{v}({\mu})$
will denote a geodesic in $X_v$ joining $\phi_v{(a)}$ and $\phi_v{(b)}$.
Let $X_{v_0} = Y$.

For convenience of exposition, we shall modify $X, {X_v}, {X_e}$ by
 quasi-isometric perturbations. Given a geodesically complete metric
space $(Z,d)$ of bounded geometry, choose a maximal disjoint collection  \\
$\{{N_1}({z_{\alpha}})\}$ of disjoint 1-balls. Then by maximality,
for all $z \in Z$ 
there exist $z_{\alpha}$ in the collection such that $d(z,{z_{\alpha}}) < 2$.
Construct a graph $Z_1$ with vertex set $\{{z_{\alpha}}\}$ and edge set consisting
of distinct vertices ${z_{\alpha}}$, ${z_{\beta}}$ 
such that $d({z_{\alpha}},{z_{\beta}}) \leq 4 $. Then $Z_1$ equipped with
the path-metric is quasi-isometric to $(Z,d)$. All metric spaces in this
section will henceforth be assumed to be graphs of edge length 1 and maps
between them will be assumed to be cellular.

We start with a  general lemma  about hyperbolic metric
spaces. This follows
easily from the fact that local quasigeodesics in a hyperbolic metric
space
are quasigeodesics \cite{GhH}.  If $x, y$ are points in a hyperbolic
metric space, $[x,y]$ will denote a geodesic joining them.

\begin{lemma}
Given $\delta > 0$,
 there exist $D, C_1$ such that  if $a, b, c, d$
are vertices of a $\delta$-hyperbolic metric space $(Z,d)$,
with ${d}(a,[b,c])={d}(a,b)$,
${d}(d,[b,c])={d}(c,d)$ and ${d}(b,c)\geq{D}$
then $[a,b]\cup{[b,c]}\cup{[c,d]}$ lies in a $C_1$-neighborhood of
any geodesic joining 
$a, d$.
\label{perps}
\end{lemma}

\smallskip

Given  a geodesic segment $\lambda\subset{Y}$, we 
now construct a quasi-convex
set
$B_{\lambda}\subset{X}$.

\medskip

{\bf Construction of $B_\lambda$}

\bigskip

Choose $C_2\geq{0}$ such that for all $e\in{T}$,
${f_e}({X_e}{\times}\{{0}\})$ and 
${f_e}({X_e}{\times}\{{1}\})$ are 
 $C_2$-quasiconvex in the appropriate vertex groups.
 Let $C{=}C_1{+}C_2$, where $C_1$ is as in Lemma \ref{perps}.

\bigskip

For $Z\subset{X_v}$, let ${N_C}(Z)$ denote the $C$-neighborhood of $Z$,
that is the set of points at distance less than or equal to $C$ from $Z$.

\smallskip
 
{\underline{\it Step 1 :}} Let 
$\mu\subset{X_v}$ be a geodesic segment in $({X_v},{d_v})$. Then
$P({\mu}) = v$. For each edge $e$ incident on $v$, but not lying on the
geodesic (in $T$) from $v_0$ to $v$, choose $p_e$, $q_e$ $\in {N_C}({\mu}){\cap}{f_v}({X_e})$ 
such that ${d_v}({p_e},{q_e})$ 
is maximal. Let ${v_1},{\cdots},{v_n}$
be terminal vertices of edges $e$ for which ${d_v}({p_e},{q_e}) > D$,
where $D$ is as in Lemma \ref{perps} above.
Observe that there are only finitely many $v_i$'s as $\mu$ is finite.
Define

\begin{center}

${B^1}({\mu}) = {i_v}({\mu}){\cup}{\bigcup}_{k=1{\cdots}n}{{\Phi}_{v_i}}({\mu})$

\end{center}

Note that $P({B^1}{({\mu})})\subset{T}$ is a  finite tree.
\medskip

{\underline{\it Step 2 :}} Step 1 above constructs 
$B^1({\lambda})$ in particular. We proceed inductively.
Suppose that $B^m({\lambda})$ has been constructed 
such that the convex hull of
$P({B^m({\lambda})}) \subset {T}$ 
is a  finite tree.
Let $ \{{w_1}, \cdots , {w_n}\} = 
P({B^m({\lambda})}){\setminus}P({B^{m-1}}({\lambda}))$. 
(Note that $n$ may depend on $m$, but we avoid repeated indices for
notational convenience.) Assume further that 
${P^{-1}}({v_k}){\cap}{B^m({\lambda})}$
is a path of the form ${i_{v_k}}({\lambda}_{k})$,
where $\lambda_k$ is a geodesic in $({X_{v_k}},{d_{v_k}})$.
Define

\begin{center}

${B^{m+1}}({\lambda})$ = 
${B^m}({\lambda}){\cup}{\bigcup}_{k=1\cdots{n}}(B^1({\lambda_k}))$

\end{center}

where $B^1({\lambda_k})$ is defined in Step 1 above.

Since each $\lambda_k$ is a finite geodesic segment in $\Gamma_H$,
the convex hull of
$P({B^{m+1}}{\lambda})$ is a  finite subtree of $T$. Further,
$P^{-1}{(v)}{\cap}B^{m+1}({\lambda})$ 
is of the form ${i_v}({\lambda_v})$ for all $v\in{P({B^{m+1}({\lambda})})}$.
This enables us to continue inductively.
Define 
\begin{center}

$B_{\lambda} = {\cup}_{m\geq{0}}B_{\lambda}^m$.

\end{center}

Note finally that the convex hull of $P({B_{\lambda}})$ in $T$ is a locally
finite tree $T_1$.

\bigskip

{\bf Quasiconvexity of $B_\lambda$}

\bigskip

We shall now show that there exists  
$C^{\prime}\geq{0}$
such that for every geodesic segment $\lambda\subset{Y}$, 
$B_{\lambda}{\subset}{X}$ is $C^{\prime}$-quasiconvex. To do this 
we construct a retraction $\Pi_{\lambda} $ from (the vertex set of)
$X$ onto ${B_\lambda}$
and show that there exists $C_0\geq{0}$ such that 
${d_X}({\Pi_\lambda}(x),{\Pi_\lambda}(y))\leq C_0{d_X}(x,y)$.
Let $\pi_v : X_v\rightarrow{\lambda_v}$ be a nearest point projection
of $X_v$ onto $\lambda_v$. $\Pi_\lambda$ is defined on
${\bigcup}_{v\in{T_1}}X_v$
by 

\begin{center}

$\Pi_{\lambda}(x) = {i_v}{\cdot}{\pi_v}(x) $ for $x\in{X_v}$.

\end{center}

If $x\in{P^{-1}}(T\setminus{T_1})$ choose $x_1\in{P^{-1}}({T_1})$
such that ${d}(x,{x_1}) = {d}(x,{P^{-1}}({T_1}))$ and define
${\Pi_{\lambda}^{\prime}}(x) = x_1$. Next define
$\Pi_\lambda{(x)} = {\Pi_\lambda}\cdot{\Pi_{\lambda}^{\prime}(x)}$.

The following Lemma says nearest point projections in a $\delta$-hyperbolic
metric space do not increase distances much.

\begin{lemma}
Let $(Y,d)$ be a $\delta$-hyperbolic metric space
 and  let $\mu\subset{Y}$ be
 a geodesic segment.
Let ${\pi}:Y\rightarrow\mu$ map $y\in{Y}$ to a point on
$\mu$ nearest to $y$. Then $d{(\pi{(x)},\pi{(y)})}\leq{C_3}d{(x,y)}$ for
all $x, y\in{Y}$ where $C_3$ depends only on $\delta$.
\label{easyprojnlemma}
\end{lemma}

{\it Proof:}
Let $[a,b]$ denote a geodesic edge-path joining vertices $a, b$. Recall that
the Gromov inner product $(a,b{)}_c$=1/2[$d{(a,c)}+{d}(b,c)-{d}(a,b)]$.
It suffices by repeated use of the triangle inequality to prove the Lemma
when ${d}(x,y)=1$. 
Let $u, v, w$ be points on $[x,\pi{(x)}]$, $[{\pi}(x),{\pi}(y)]$ and 
$[{\pi}(y),x]$ respectively such that ${d}(u,{\pi}(x))={d}(v,{\pi}(x))$,
${d}(v,{\pi}(y))={d}(w,{\pi}(y))$ and ${d}(w,x)={d}(u,x)$.
Then ${(x,{\pi}(y))}_{\pi{(x)}}={d}(u,{\pi}(x))$. Also, since $Y$
 is $\delta$-hyperbolic, 
the diameter of the inscribed triangle with vertices $u, v, w$ is
less than or equal to $2\delta$ (See \cite{Shortetal}).

\begin{eqnarray*}
{d}(u,x)+{d}(u,v) & \geq & {d}(x,{\pi}(x)) =  {d}(u,x)+{d}(u,{\pi}(x)) \\
\Rightarrow {d}(u,{\pi}(x)) & \leq & {d}(u,v)\leq{2\delta}  \\
\Rightarrow {(x,{\pi}(y))}_{{\pi}(x)} & \leq & 2{\delta}
\end{eqnarray*}

Similarly, ${(y,{\pi}(x))}_{\pi{(y)}}\leq{2\delta}$.

\begin{center}
i.e. ${d}(x,{\pi}(x))+{d}({\pi}(x),{\pi}(y))-{d}(x,{\pi}(y))\leq{4\delta}$

and  ${d}(y,{\pi}(y))+{d}({\pi}(x),{\pi}(y))-{d}(y,{\pi}(x))\leq{4\delta}$
\end{center}

Therefore,
\begin{eqnarray*}
\lefteqn{2{d}({\pi}(x),{\pi}(y)) }   \\
    & \leq & {8\delta}+{d}(x,{\pi}(y))-{d}(y,{\pi}(y))+{d}(y,{\pi}(x))-{d}(x,{\pi}(x)) \\
    & \leq & {8\delta}+{d}(x,y)+{d}(x,y) \\
     & \leq & {8\delta}+2     
\end{eqnarray*}
 Hence  ${d}({\pi}(x),{\pi}(y))\leq{4\delta}+1$.
Choosing $C_3 = {4\delta}+1$, we are through. $\Box$

\medskip

\begin{lemma}
Let $(Y,d)$ be a $\delta$-hyperbolic metric space.
Let $\mu$ be a geodesic segment in $Y$ with end-points $a, b$ and let
$x$ be any vertex in $Y$. Let $y$ be a vertex on $\mu$ such that 
${d}(x,y)\leq{d}(x,z)$ for any $z\in\mu$. Then a geodesic path
from $x$ to $y$ followed by a geodesic path from $y$ to $z$ is a 
$k$-quasigeodesic for some $k$ dependent only on $\delta$.
\label{unionofgeodslemma}
\end{lemma}
{\it Proof:}
As in Lemma \ref{easyprojnlemma} let $u, v, w$ be points on edges $[x,y]$,
$[y,z]$ and $[z,x]$ respectively such that ${d}(u,y)={d}(v,y)$,
 ${d}(v,z)={d}(w,z)$ and  ${d}(w,x)={d}(u,x)$.
Then ${d}(u,y)={(z,x)}_y\leq{2\delta}$ and the inscribed triangle with vertices
$u, v, w$ has diameter less than or equal to $2\delta$ (See \cite{Shortetal}).
 $[x,y]\cup{[y,z]}$ is a union of 2 geodesic paths lying in a 
$4\delta$ neighborhood of a geodesic $[x,z]$. Hence a geodesic path from
$x$ to $y$ followed by a geodesic path from $y$ to $z$ is a $k-$quasigeodesic
for some $k$ dependent only on $\delta$. $\Box$
\medskip
\begin{lemma}
Suppose $(Y,d)$ is a $\delta$-hyperbolic metric space.
If $\mu$ is a $({k_0},{\epsilon_0})$-quasigeodesic in $Y$ and $p, q, r$ are 
3 points in order on $\mu$ then ${(p,r)}_q\leq{k_1}$ for some $k_1$
dependent on $k_0$, $\epsilon_0$ and $\delta$ only.
\label{qgeodiplemma}
\end{lemma}

{\it Proof:}  $[a,b]$ will denote a geodesic path joining $a, b$. 
Since $p, q, r$ are 3 points in order on $\mu$, $[p,q]$ followed
by $[q,r]$ is a $({k_0},{\epsilon_0})$-quasigeodesic in the $\delta$-hyperbolic metric space $Y$. 
Hence there exists a $k_1$ dependent on $k_0$, $\epsilon_0$
 and $\delta$ alone
such that ${d}(q,[p,r])\leq{k_1}$. Let $s$ be a point on $[p,r]$ such that
${d}(q,s)={d}(q,[p,r])\leq{k_1}$. Then 
\begin{eqnarray*}
{(p,r)}_q & = & 1/2({d}(p,q)+{d}(r,q)-{d}(p,r))  \\
           & = & 1/2({d}(p,q)+{d}(r,q)-{d}(p,s)-{d}(r,s)) \\
         & \leq & {d}(q,s)\leq{k_1}. \Box
\end{eqnarray*}

\medskip

\medskip

\begin{lemma}
Suppose $(Y,d)$ is $\delta$-hyperbolic.
Let $\mu_1$ be some geodesic segment in $Y$ joining $a, b$ and let $p$
be any vertex of $Y$. Also let $q$ be a vertex on $\mu_1$ such that
${d}(p,q)\leq{d}(p,x)$ for $x\in\mu_1$. 
Let $\phi$ be a $(K,{\epsilon})$ - quasiisometry from $Y$ to itself.
Let $\mu_2$ be a geodesic segment 
in $Y$ joining ${\phi}(a)$ to ${\phi}(b)$ for some $g\in{S}$. Let
$r$ be a point on $\mu_2$ such that ${d}({\phi}(p),r)\leq{d}({\phi(p)},x)$ for $x\in\mu_2$.
Then ${d}(r,{\phi}(q))\leq{C_4}$ for some constant $C_4$ dependent   only on
$K, \epsilon $ and $\delta$. 
\label{cruciallemma}
\end{lemma}

{\it Proof:}
Since  ${\phi}({\mu_1})$ is a 
$({K,\epsilon})-$   quasigeodesic
joining ${\phi}(a)$ to ${\phi}(b)$, it lies in a $K^{\prime}$-neighborhood 
of $\mu_2$ where $K^{\prime}$ depends only on $K, {\epsilon}, \delta$. 
Let $u$ be a vertex
in ${\phi}({\mu_1})$ lying at a distance at most $K^{\prime}$ from $r$.
Without loss of generality suppose that $u$ lies on ${\phi}([q,b])$, 
where $[q,b]$ denotes the geodesic subsegment of $\mu_1$ joining $q, b$.
[See Figure 1.]

\begin{figure}
\hbox to \hsize{\hss\psfig{file=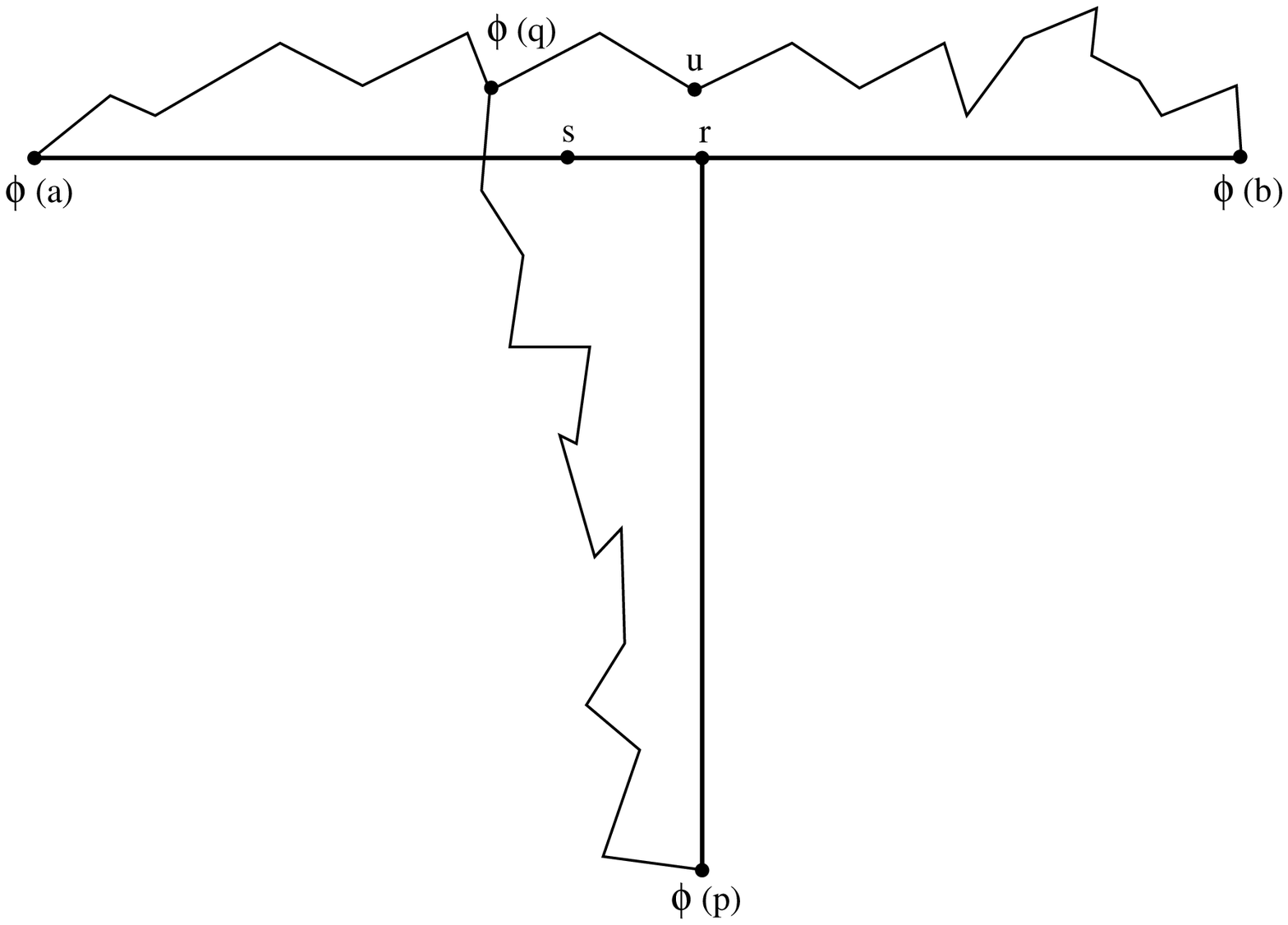,width=5in}\hss}
\caption{}
\end{figure}

Let $[p,q]$ denote a geodesic joining $p, q$.
From Lemma \ref{unionofgeodslemma} $[p,q]\cup{[q,b]}$ is a $k$-quasigeodesic,
where $k$ depends on $\delta$ alone. Therefore 
${\phi}([p,q])\cup{\phi}([q,b])$ is a 
$({{K_0}, \epsilon_0})$-quasigeodesic, where $K_0, {\epsilon_0}$ depend
on $K, k, \epsilon$. Hence, by Lemma \ref{qgeodiplemma} 
${({\phi}(p),u)}_{{\phi}(q)}\leq{K_1}$, where $K_1$ depends on $K$, $k$,
$\epsilon$ and $\delta$ alone. Therefore, 
\begin{eqnarray*}
\lefteqn{ {({\phi}(p),r)}_{{\phi}(q)} } \\
  & = & 1/2[{d}({\phi}(p),{\phi}(q))+{d}(r,{\phi}(q))-{d}(r,{\phi}(p))]  \\
  & \leq & 1/2[{d}({\phi}(p),{\phi}(q))+{d}(u,{\phi}(q))+{d}(r,u) \\ 
  &      &  \hspace{1.5in}       -{d}(u,{\phi}(p))+{d}(r,u)] \\
  & = & {({\phi}(p),u)}_{{\phi}(q)}+{d}(r,u) \\
  & \leq & {K_1}+{K^{\prime}}
\end{eqnarray*}

There exists $s\in{\mu_2}$ such that ${d}(s,{\phi}(q))\leq{K^{\prime}}$

\begin{eqnarray*}
{{({\phi}(p),r)}_{s}} & = & 1/2[{d}({\phi}(p),s)+{d}(r,s)-{d}(r,{\phi}(p))] \\
   & \leq &  1/2[{d}({\phi}(p),{\phi}(q))+{d}(r,{\phi}(q))-{d}(r,{\phi}(p))]+{K^{\prime}}  \\
  & = & {({\phi}(p),r)}_{{\phi}(q)}+{K^{\prime}} \\
  & \leq &  {K_1}+{K^{\prime}}+{K^{\prime}} \\
  & = &  {K_1}+2{K^{\prime}}
\end{eqnarray*}

Also, as in the proof of Lemma \ref{easyprojnlemma} 
${({\phi}(p),s)}_r\leq{2\delta}$

\begin{eqnarray*}
{d}(r,s) & = & {({\phi}(p),s)}_{r}{+}{({\phi}(p),r)}_{s}  \\
             & \leq & K_1{+}2{K^{\prime}}{+}2{\delta} \\
{d}(r,{\phi}(q)) & \leq & K_1{+}2{K^{\prime}}{+}2{\delta}{+}{d}(s,{\phi}(q)) \\
                         & \leq & K_1{+}2{K}^{\prime}{+}2{\delta}{+}{K}^{\prime}  
\end{eqnarray*}

Let $C_4{=}K_1{+}3{K^{\prime}}{+}2{\delta}$. Then ${d}(r,{\phi}(q))\leq{C_4}$ 
and $C_4$ is independent of $a, b, p$. $\Box$

\medskip

$d_T$ will denote the metric on $T$. We are now in a position to prove:

\begin{theorem}
There exists $C_0\geq{0}$ such that 
${d}({\Pi_\lambda}(x),{\Pi_\lambda}(y))\leq C_0{d_G}(x,y)$ for $x, y$
vertices of $\Gamma_G$.
\label{mainref}
\end{theorem}

{\it Proof:} It suffices to prove 
the theorem when ${d_G}(x,y)=1$.  

\medskip

{\it Case (a):} $x, y \in {P^{-1}}(v)$ for some $v\in{T_1}$. 
 From Lemma \ref{easyprojnlemma}, there exists $C_3$ such that
${d_v}({\pi_v}{\cdot}{i_v^{-1}}(x),{\pi_v}{\cdot}{i_v^{-1}}(y))\leq{C_3}$.
Since embeddings of $X_v$ in $X$ are cellular, 
$d({\Pi_{\lambda}}(x),{\Pi_{\lambda}}(y)) \leq C_3$.
	
{\it Case (b):} $x\in{P^{-1}}(w)$ and $y\in{P^{-1}}(v)$
for some $v, w\in{T_1}$.

Since ${d_G}(x,y) = 1$, $v, w$ are adjacent in $T_1$.  Assume, without
loss of generality, $w = v_{-}$.

Recall that 

\begin{center}

${B_{\lambda}}\cap{P^{-1}}(v) = {i_v}({\lambda_v})$ \\
${B_{\lambda}}\cap{P^{-1}}(w) = {i_w}({\lambda_w})$ \\

\end{center}

Also,  $\lambda_v$ = ${\Phi_v}({\mu_w})$, for some geodesic $\mu_w$ contained
in $X_w$, such that end-points of $\mu_w$ lie in a $C$-neighborhood of
$\lambda_w$.

Let $z\in{(X_w)}$ denote a nearest point projection of ${i_w^{-1}}(x)$
 onto $\mu_w$. 
Then, by Lemma \ref{cruciallemma}, 
\begin{center}
$d({i_w}(z),{{\Pi_{\lambda}}}\cdot{\phi_v}(x)) \leq  
d({i_w}(z),{\phi_v}{\cdot}{i_w}(z)) 
+ d({\phi_v}{\cdot}{i_w}(z),{{\Pi_{\lambda}}}\cdot{\phi_v}(x))
 \leq 1 + C_4$.
\end{center}

Since, $d(x,y) = 1 = d(x,{\phi_v}(x))$ and $i_v$'s are uniformly proper
embeddings, there exists $C_5 > 0$ such that ${d_v}({\phi_v}(x),y) \leq C_5$
and \\
 $d({{\Pi_{\lambda}}}({\phi_v}(x)),{{\Pi_{\lambda}}}(y)) \leq {C_3}{C_5}$. 

Since the end-points of $\mu_w$ lie in a $C$-neighborhood of $\lambda_w$,
there exists $C_6$, depending on $\delta$ and $C$ such that
 $d(z,{{\Pi_{\lambda}}}(x)) \leq C_6$.

Finally, by the triangle inequality, 
\begin{center}
$d({{\Pi_{\lambda}}}(x),{{\Pi_{\lambda}}}(y)) \leq C_6 + 1 + C_4 + {C_3}{C_5} = {C_7}(say)$
\end{center}

\medskip

{\it Case (c):} $P([x,y])$ is not contained in $T_1$. 

Since ${d}(x,y) = 1$
$P(x)$ and $P(y)$ belong to the closure $T_2$ of the 
same component  of $T{\setminus}T_1$.
Then ${P}{\cdot}{{\Pi}_{\lambda}^{\prime}}(x) =
{P}{\cdot}{{\Pi}_{\lambda}^{\prime}}(y) = v$ for some $v\in{T}$.

Also ${d}({\Pi_{\lambda}}(x),{\Pi_{\lambda}}(y)) = 
{d}({\Pi_{\lambda}}{\cdot}{\Pi_{\lambda}^{\prime}}(x),{\Pi_{\lambda}}{\cdot}{\Pi_{\lambda}^{\prime}}(y)) $

Let ${x_1} = \Pi_{\lambda}^{\prime}(x) $ and
${y_1} = \Pi_{\lambda}^{\prime}(y) $.

Let $D$ and $C_1$
 be as in Lemma \ref{perps}. If ${d}({\Pi}_{\lambda}({x_1}),{\Pi}_{\lambda}({y_1})) \geq D$, let 

\begin{eqnarray*}
{i_v}^{-1}({x_1}) & = & u_1, \\
{i_v}^{-1}({\Pi_{\lambda}}({x_1})) & = & u_2, \\
{i_v}^{-1}({y_1}) & = & v_1, \\
{i_v}^{-1}({\Pi_{\lambda}}({y_1})) & = & v_2.
\end{eqnarray*}

Then by 
Lemma \ref{perps} 
$[{u_1},{u_2}]{\cup}[{u_2},{v_2}]{\cup}[{v_2},{v_1}]$
is a quasigeodesic lying in a $C_1$-neighborhood of $[{u_1},{v_1}]$.

Also, ${x_1}, {y_1}\in{i_v}({X_v})$.
Since the image of an edge space in a vertex space
 is $C_2$-quasiconvex,  
there exist $e\in{T}$ and ${x_2}, {y_2}{\in}{f_e}({X_e}{\times}\{{0}\})$
 such that 
${d}({x_2},{\Pi_{\lambda}}({x_1})) \leq C_1 + C_2 = C$
and ${d}({y_2},{\Pi_{\lambda}}({y_1})) \leq C_1 + C_2 = C$.

By construction 
${d}({\Pi_{\lambda}}({x_2}),{\Pi_{\lambda}}({y_2})) \leq D$.
(Else the edge $P([x,y])$ of $T$ would be in $T_1$).
Therefore,

\begin{eqnarray*}
{d}({\Pi_{\lambda}}(x), {\Pi_{\lambda}}(y)) & =   &
{d}({\Pi_{\lambda}}({x_1}), {\Pi_{\lambda}}({y_1})) \\
   & \leq & 2{C} + D + 2{C} \\
    &  =  &  4{C} + D
\end{eqnarray*}

Choosing ${C_0} =$ max $\{{C_3}, {C_{7}}, {4{C}} + D\}$, we are through.
$\Box$

\medskip

To complete the proof of our main Theorem, we need a final Lemma.

\begin{lemma}
There exists  $A > 0$, such that 
if $a\in{P^{-1}}(v){\cap}{B_{\lambda}}$ for some $v\in{T_1}$ then
there exists $b\in{i({\lambda})}$ with
${d}(a,b)$ $ \leq $ ${A}{d_T}(Pa,Pb)$.
\label{connectionlemma}
\end{lemma}
 
{\it Proof:}
Let $\mu$ be a geodesic path from $v_0$ to $v$ in $T$.
 Order the vertices on $\mu$ so that we have a finite
sequence ${v_0}={y_0},{y_1},{\cdots},{y_n}=v$ such that ${d_T}({y_i},{y_{i+1}})=1$.
and ${d_T}({v_0},v)=n$. Recall further, $P({B_{\lambda}}) = T_1$. Hence
$y_i \in T_1$.

Recall that $B_{\lambda}$ is of the form 
${\bigcup}_{v\in{T_1}}{i_v}({\lambda_v})$.

It suffices to prove that there exists $A > 0$ independent of $v$
such that if 
$p \in {{i_{y_j}}({\lambda_{y_j}})}$, there exists 
$q \in {i_{y_{j-1}}}({\lambda_{y_{j-1}}})$
with $d(p,q) \leq A$.

By construction, $\lambda_{y_j} = {\Phi_{y_j}}({\mu})$ for some geodesic
$\mu $ in $X_{y_{j-1}}$ such that end-points of $\mu$ lie in a 
$C$-neighborhood of $\lambda_{y_{j-1}}$. 
Since $\phi_{y_j}$ is a quasi-isometry,
there exists $C_1$ such that $p$ lies in a $C_1$ neighborhood of 
$\phi_{y_j}{(q_0)}$ for some $q_0{\in}\mu$. Therefore, $d({q_0},p) \leq 1 + C$.

Also, since 
 end-points of $\mu$ lie in a 
$C$-neighborhood of $\lambda_{y_{j-1}}$, there exists 
$q \in {i_{y_j}}({\lambda_{y_{j-1}}})$
with $d({q_0},q) \leq C_2$ where $C_2$ depends only on $\delta$ and $C$.
Choosing $A = 1 + C + C_2$, we are through. $\Box$
\medskip

The main theorem of this paper follows:

\begin{theorem}
 Let (X,d) be a tree (T) of hyperbolic metric spaces satisfying the
qi-embedded condition.  Let $v$ be a vertex of $T$. Then 
${i_v} : X_v \rightarrow X$ extends continuously to ${\hat{{i_v}}} : \widehat{X_v}
\rightarrow \widehat{X}$.
\label{main}
\end{theorem}

{\it Proof:}  Without loss of generality, let $v_0 = v$ be the base vertex of $T$.
To prove the existence of a Cannon-Thurston map,
it suffices to show  (from  Lemma \ref{contlemma})
that for all $M\geq{0}$ and ${x_0} \in X_v$ there exists $N\geq{0}$
such that if a geodesic segment 
$\lambda$ lies oustside the $N$-ball around ${x_0}\in{X_v}$,
$B_{\lambda}$ lies outside the $M$-ball around ${i_v}({x_0})\in{X}$.

To prove this, we show that if $\lambda$ lies outside the 
$N$-ball around ${x_0}\in{X_v}$,
$B_{\lambda}$ lies outside a certain  $M(N)$-ball around ${i_v}({x_0})\in{X}$,
 where
$M(N)$ is a proper function from $\Bbb{N}$ into itself.

Since $X_v$ is properly embedded in $X$ there exists $f(N)$
such that ${i_v}({\lambda})$ lies outside the $f(N)-$ball around $x_0$ 
in $X$
and $f(N)\rightarrow\infty$ as $N\rightarrow\infty$.

Let $p$ be any point on $B_{\lambda}$. There exists $y\in{{i_v}({\lambda})}$
such that ${d}(y,p)\leq{A}{d_T}(Py,Pp)$ by Lemma \ref{connectionlemma}.
Therefore,
\begin{eqnarray*}
{d}({x_0},p) & \geq & {d}({x_0},y)-A{d_T}(Py,Pp) \\
                     & \geq & f(N)-A{d_T}(P({x_0}),Pp)
\end{eqnarray*}

By our choice of metric on $X$,

\begin{center}
${d}({x_0},p) \geq {d_T}(P({x_0}),Pp) $
\end{center}

Hence
\begin{eqnarray*}
 {d}({x_0},p) & \geq & max({f(N)-A{d_T}(P({x_0}),Pp)}, {d_T}(P({x_0}),Pp)) \\
                 & \geq & \frac{f(N)}{A+1}
\end{eqnarray*}

>From Theorem \ref{mainref}
there exists $C^{\prime}$ independent of $\lambda$ such that
 $B_{\lambda}$ is a $C^{\prime}$-quasiconvex set containing ${i_v}({\lambda})$.
Therefore
any geodesic joining the end-points of ${i_v}({\lambda})$ lies in a 
$C^{\prime}$-neighborhood of $B_{\lambda}$.

Hence any geodesic joining  end-points of ${i_v}({\lambda})$ lies outside a ball of radius $M(N)$ where
\begin{center}
$M(N) = {\frac{f(N)}{A+1}}-{C^{\prime}}$
\end{center}
Since $f(N)\rightarrow\infty$ as $N\rightarrow\infty$ so does $M(N)$.
$\Box$

\medskip

The following is a direct consequence of Theorem \ref{main} above.

\begin{cor}
 Let $G$ be a hyperbolic group acting cocompactly on a simplicial tree
$T$ such that all vertex and edge stabilizers are hyperbolic. Also
suppose that every inclusion of an edge stabilizer in a vertex stabilizer
is a quasi-isometric embedding. Let $H$ be the stabilizer of a vertex or
edge of $T$. Then 
there exists a Cannon-Thurston map for $(H,G)$.
\label{main1}
\end{cor}

\section{Geometrically Tame Kleinian Groups}

In this section we apply Theorem \ref{main} to geometrically tame Kleinian 
groups.

The {\it convex core} of a hyperbolic 3-manifold $N$ (without cusps)
is the smallest convex submanifold $C(N) \subset N$ for which inclusion 
is a homotopy equivalence. If $C(N)$ has finite volume, $N$ is said to be
{\it geometrically finite}. There exists a compact 3-dimensional submanifold
$M \subset N$, the {\it Scott core} \cite{Scott1} whose inclusion is
a homotopy equivalence. 
The ends of $N$ are in one-to-one correspondence with the components of
$N - M$ or, equivalently, the components of $\partial{M}$. We say that an
end of $N$ is {\it geometrically finite} if it has a neighborhood missing
$C(N)$. An end of $N$ is {\it simply degenerate} if it has a neighborhood
homeomorphic to $S{\times}{\Bbb{R}}$, where $S$ is the corresponding
component of $\partial{M}$, and if there is a sequence of pleated surfaces
homotopic in this neighborhood to the inclusion of $S$, and exiting every
compact set. $N$ is called {\it geometrically tame} if all of its ends are
either geometrically finite or simply degenerate. In particular, $N$ is homeomorphic
to the interior of $M$. For a more detailed discussion of pleated surfaces 
and geometrically tame ends, see \cite{Thurstonnotes} or \cite{Minsky2}.

Let ${inj}_N(x)$ denote the injectivity radius at $x\in{N}$. For the purposes
of this section, we shall assume that there exists $\epsilon_0 > 0$ such that
${inj}_N{(x)} > \epsilon_0$ for all $x\in{N}$. In order to apply Theorem \ref{main}
we need some preliminary Lemmas.

Let $E$ be a simply degenerate end of $N$. Then $E$ is homeomorphic to
$S{\times}[0,{\infty})$ for some closed surface $S$ of genus greater than one.

\begin{lemma}
\cite{Thurstonnotes}
There exists $D_1 > 0$ such that for all $x\in{N}$, there exists a pleated 
surface $g : (S,{\sigma}) \rightarrow N$ with $g(S){\cap}{B_{D_1}}(x) \neq \emptyset$.
\label{closepleated}
\end{lemma}

The following Lemma follows easily from the fact that ${inj}_N{(x)} > \epsilon_0$:

\begin{lemma} \cite{Bonahon1},\cite{Thurstonnotes}
There exists $D_2 > 0$ such that if $g : (S,{\sigma}) \rightarrow N$ is a
pleated surface, then $dia(g(S)) < D_2$. 
\label{diameter}
\end{lemma}

The following Lemma due to Minsky \cite{Minsky2} follows from compactness
of pleated surfaces.
\begin{lemma}
\cite{Minsky2}
Fix $S$ and $\epsilon > 0$. Given $a > 0$ there exists $b > 0$ such that if
$g : (S,{\sigma})\rightarrow{N}$
and $h : (S,{\rho})\rightarrow{N}$ are homotopic pleated surfaces which
are isomorphisms on $\pi_1$ and ${inj}_N{(x)} > \epsilon$ for all $x \in N$,
then\\
\begin{center}
${d_N}(g(S),h(S)) \leq a \Rightarrow {d_{Teich}}({\sigma},{\rho}) \leq b$,
\end{center}
    where $d_{Teich}$ denotes Teichmuller distance.
\label{pleatedcpt}
\end{lemma}

\begin{lemma}
There exist $K, \epsilon$ and  a homeomorphism $h$ from $E$ to the
universal curve over a Lipschitz path in Teichmuller space, such that $h$
is a $(K,{\epsilon})$-quasi-isometry.
\label{lipcurve}
\end{lemma}

{\bf Proof:}
We can assume that $S{\times}\{{0}\}$ is mapped to a pleated surface
$S_0$ $\subset N$ under the homeomorphism from $S{\times}[0,{\infty})$
to $E$. We shall construct inductively a sequence of `equispaced' pleated
surfaces ${S_i}\subset{E}$ exiting the end. Assume that ${S_0},{\cdots},{S_n}$
have been constructed such that:
\begin{enumerate}
\item If $E_i$ be the non-compact component of $E{\setminus}{S_i}$, then
$S_{i+1} \subset E_i$.
\item Hausdorff distance between $S_i$ and $S_{i+1}$ is bounded above by
$3({D_1}+{D_2})$.
\item ${d_N}({S_i},{S_{i+1}}) \geq D_1 + D_2$.
\item From Lemma \ref{pleatedcpt} and condition (2)
above there exists $D_3$ depending on $D_1$, $D_2$ and $S$ such that
$d_{Teich}({S_i},{S_{i+1}}) \leq D_3$
\end{enumerate}

Next choose $x \in E_n$, such that ${d_N}(x,{S_n}) = 2({D_1}+{D_2})$. Then
by Lemma \ref{closepleated}, there exists a pleated surface
$g : (S,{\tau}) \rightarrow N$ such that 
${d_N}(x,{g(S)}) \leq D_1$. Let ${S_{n+1}} = g(S)$. Then by the triangle
inequality and Lemma \ref{diameter}, if $p\in{S_n}$ and $q\in{S_{n+1}}$,
\begin{center}
${D_1} + {D_2} \leq {d_N}(p,q) \leq 3({D_1} + {D_2})$.
\end{center}

This allows us to continue inductively. The Lemma follows. $\Box$

Observe that $\widetilde{E}$ is quasi-isometric to a tree (in fact a ray)
of hyperbolic metric spaces by setting $T = [0,{\infty})$, with vertex
set $\{n : n \in {\Bbb{N}}{\cup}\{{0}\} \}$, edge set 
$\{ [{n-1},n]: n\in{\Bbb{N}} \}$,
${X_n} = {\widetilde{S_n}} = X_{[n-1,n]}$.  Further, by Lemma \ref{pleatedcpt}
this tree of hyperbolic metric spaces satisfies the quasi-isometrically
embedded condition.  We shall now describe $\widetilde{C(N)}$ as a tree of
hyperbolic metric spaces. Assume $M{\subset}C(N)$ and $\partial{M} = \{{F_1},{\cdots},{F_n}\}$
where $F_i$ are pleated surfaces in $N$ cutting off ends $E_i$.

\begin{lemma}
\cite{BF}
$\pi_1{(N)}$ is hyperbolic in the sense of Gromov. Also, if $i : E \rightarrow
N$, denotes inclusion, then ${i_*}{\pi_1}(E)$ is a quasiconvex subgroup
of ${\pi_1}(N)$.
\label{hyppi1}
\end{lemma}

{\bf Remark:} In fact there exists a geometrically finite hyperbolic manifold
homeomorphic to $N$. This is part of Thurston's monster theorem. See 
\cite{ctm} for a different proof of the fact. Also, the limit set of
a geometrically finite manifold is locally connected \cite{and-mask}.
This shall be of use later.

\medskip

Note that $\widetilde{M} \subset \widetilde{N}$ is quasi-isometric
 to the Cayley graph of $\pi_1{(N) }$. Hence, $\widetilde{M}$ is a hyperbolic
metric space. Let $\widetilde{F_i} \subset \widetilde{N}$ represent a lift
of $F_i$ to $\widetilde{N}$. Then, by Lemma \ref{hyppi1} above, 
$\widetilde{F_i}$ is a word-hyperbolic metric space. If $\widetilde{E_i}$ is
a lift of $E_i$ containing $\widetilde{F_i}$ then from our previous discussion,
$\widetilde{E_i}$ is a ray of hyperbolic metric spaces. Since there are only
finitely many ends $E_i$, we have thus shown:

\begin{lemma}
The hyperbolic metric space $\widetilde{C(N)}$ is quasi-isometric to a tree
(T) of hyperbolic metric spaces satisfying the qi-embedded condition. Further,
we can choose a base vertex $v_0$ of $T$ such that $X_{v_0}$ is homeomorphic
to $\widetilde{M}$.
\label{hyptree}
\end{lemma}

Applying Theorem \ref{main}, we get 

\begin{theorem}
 Let $\Gamma$ be a geometrically tame Kleinian 
group, such that ${{\Bbb{H}}^3}/{\Gamma} = M$ has injectivity radius 
uniformly bounded below by some $\epsilon > 0$. Then there exists a continuous
map from the Gromov boundary of $\Gamma$ (regarded as an abstract group)
to the limit set of $\Gamma$ in ${\Bbb{S}}^2_{\infty}$.
\label{main2}
\end{theorem}

The above theorem has been independently proven by Klarreich \cite{klarreich}
using different techniques.

\begin{lemma} 
Let $N$ be a geometrically tame 3-manifold with 
${inj}_N{(x)} > \epsilon_0 > 0$ for all $x \in N$. Then the Gromov boundary of 
$\pi_1{(N)}$ is locally connected.
\label{locconnlemma}
\end{lemma}

{\bf Proof:} This follows from the fact that there exists a geometrically
finite manifold $M = {{\Bbb{H}}^3}/{\Gamma}$ 
homeomorphic to $N$ \cite{ctm} and that for such an
$M$, the limit set of  $\Gamma$ is
locally connected \cite{and-mask}. $\Box$

Since a continuous image of a compact locally connected set is locally 
connected, Lemma \ref{locconnlemma}  and Theorem \ref{main2} give:

\begin{cor}
Let $N = {{\Bbb{H}}^3}/{\Gamma}$ be a geometrically tame 3-manifold with 
${inj}_N{(x)} > \epsilon_0$ for all $x \in N$. Then the limit set of 
$\Gamma$ is locally connected.
\label{locconn}
\end{cor}

Lemma \ref{lipcurve} shows that there exists a quasi-isometry from a lift  
$\widetilde{E}$ of an end to the universal cover of a universal curve over a Lipschitz
path $\sigma$ in $Teich(S)$. We show further that $\sigma$ is a Teichmuller
quasigeodesic.

It is well known that geodesics in hyperbolic metric spaces diverge 
exponentially. The following proposition `quasi-fies' this statement:

\begin{prop}
Given ${\delta}$,  ${A_0}\geq{0}$ there exist $\beta {>} {1}$,  $B > 0$,
$K \geq 1$ and $\epsilon \geq 0$
such that if
$[x,y], [y,z]$ and $[z,w]$ are geodesics in a $\delta$-hyperbolic metric space
$(X,d)$  with 
$(x,z)_y\leq{A_0}$, $(y,w)_z\leq{A_0}$
and ${d}(y,z) \geq B$
 then   any path joining $x$ to $w$ and lying outside a 
$D$-neighborhood of $[y,z]$ has length greater than
or equal to ${\beta}^{D}{d}(y,z)$,\\ 
where $D = min{\{{({d}(x,[y,z])-1)},{({d}(w,[y,z])-1)}\}}$.\\

\label{ineffprop}
\end{prop}

\begin{lemma}
$\sigma$ is a Teichmuller quasigeodesic.
\label{teichqgeod}
\end{lemma}

{\bf Proof:}
Let $S_0 = {\partial}E$ be a pleated surface containing a closed geodesic $l$
of $N$. This can always be arranged by taking a simple closed geodesic
sufficiently far out in $E$ and mapping in a pleated surface containing
it \cite{Thurstonnotes}. Construct a sequence of equispaced pleated
surfaces as in Lemma \ref{lipcurve}. 
$\widetilde{E}$ is quasi-isometric to a ray of
hyperbolic metric spaces $X$,
with vertex
set $\{n : n \in {\Bbb{N}}{\cup}\{{0}\} \}$, edge set 
$\{ [{n-1},n]: n\in{\Bbb{N}} \}$,
${X_n} = {\widetilde{S_n}} = X_{[n-1,n]}$. 

 Fix $x_0$ in $S_0$. Inductively, define $x_n$
to be the image of $x_{n-1}$ under the Teichmuller map from $S_{n-1}$ to
$S_n$. Let $r$ denote a quasi-isometric embedding of $[0,{\infty})$
sending $[n,n+1]$ to the shortest geodesic from $x_n$ to $x_{n+1}$.
Then $r$ is a quasigeodesic in $E$. Let $[a,b]$ be a lift of $l$ to
$\widetilde{E}$. Let $\lambda \subset X$ be the image of $[a,b]$
under a quasi-isometric homeomorphism $h$ between $\widetilde{E}$ and $X$,
sending $\widetilde{S_n}$ to $X_n$.
Construct $B_{\lambda}\subset{X}$, as in the previous section.
 Lifts of $r$ through $a, b$
 diverge exponentially. 

>From 
Theorem \ref{mainref}, $B_{\lambda}$ is quasiconvex
and hence a hyperbolic metric space with the inherited metric.
 Let $r_1$, $r_2$ be the images of these
lifts through the end-points of $\lambda$. Then $r_1$, $r_2$ are 
$(K,{\epsilon})$-quasigeodesics  diverging exponentially in $X$.
Assume, after reparametrization if necessary, ${r_1}(n)$, ${r_2}(n) \in X_n$.
Let $d_n$ denote the path metric on $X_n$. Then by Proposition \ref{ineffprop},
there exist $C_1 > 1, n \geq 1$ such that 
${d_{N+n}}({r_1}(N+n),{r_2}(N+n)) \geq {C_1}{d_N}({r_1}(N),{r_2}(N))$ for
all $N \geq 0$. Hence there exists $C_2$ such that $d_{Teich}({S_N},{S_{N+n}})
\geq C_2$ for all $N \geq 0$. 
Since $\sigma$ was shown to be Lipschitz in Lemma \ref{lipcurve},
this proves that  $\sigma$ is a Teichmuller quasigeodesic. $\Box$

So far arguments have been coarse. At this stage, we need to quote 
a part of the main theorem of \cite{Minsky2}.

\begin{theorem} \cite{Minsky2}
If $N$ is a geometrically tame hyperbolic 3-manifold with indecomposable
fundamental group, such that there exists $\epsilon_0 > 0$ with
${inj}_N{(x)} > \epsilon_0$ for all $x\in{N}$, then each
simply degenerate end $E$ of $N$ gives rise to a unique Teichmuller ray $r$,
such that every pleated surface in $E$ lies at a uniformly bounded distance
from $r$. Further, $r$ depends only on the corresponding ending lamination.
\label{minskymain}
\end{theorem}

That $r$ depends only on the corresponding ending lamination was proven
by Masur \cite{masur}.

Combining Lemma \ref{teichqgeod} and Theorem \ref{minskymain} we have 
a new proof of the main theorem of \cite{Minsky} :
 the ending lamination
theorem for 3-manifolds with freely indecomposable fundamental group and 
a uniform lower bound on injectivity radius.

\begin{theorem}
Let $N_1$ and $N_2$ be homeomorphic
hyperbolic 3-manifolds with freely indecomposable fundamental group. Suppose
there exists a uniform
lower bound $\epsilon > 0$  on the injectivity radii of $N_1$ and $N_2$.
If the end invariants of corresponding ends of $N_1$ and $N_2$ are
equal, then $N_1$ and $N_2$ are isometric.
\label{main3}
\end{theorem}

{\bf Proof:} From Lemma \ref{teichqgeod}, corresponding simply degenerate
ends $E_{i1}$, $E_{i2}$ of $N_1$ and $N_2$ are homeomorphic via 
quasi-isometries  to universal curves over
Teichmuller  quasi-geodesics $l_{i1}$ and
$l_{i2}$ lying in  bounded neighbourhoods
 of Teichmuller  geodesics $l_i$. Hence corresponding
ends are homeomorphic via quasi-isometries  to each other. Therefore 
$N_1$, $N_2$ are homeomorphic by a quasi-isometry. Finally, by
\cite{Sullivan} $N_1$ and $N_2$  are isometric.

\section{Examples}

Let $H$ be a hyperbolic subgroup of a hyperbolic group $G$.

{\bf Definition :} \cite{Gromov2} \cite{farb}
{\it If $i: \Gamma_H \rightarrow \Gamma_G$ 
be an embedding of the Cayley
graph of $H$ into that of $G$, then the {\it distortion} function is
given by
\begin{center}
$disto(R) = R^{-1} Diam_{\Gamma_0}({\Gamma_0}{\cap}B(R))$,
\end{center}
where $B(R)$ is the ball of radius $R$ around $1\in{\Gamma_G}$.}

\medskip

All previously known examples of non-quasiconvex 
hyperbolic subgroups of
hyperbolic groups exhibit
exponential distortion. We construct in this section some examples
exhibiting greater distortion. Some of these will be shown to have
Cannon-Thurston maps. For the rest, existence of Cannon-Thurston maps is
not yet known. Further, we shall describe certain examples of free
subgroups of $PSL_2{\Bbb{C}}$ and show that they exhibit arbitrarily
large distortion. The existence of Cannon-Thurston maps for some of these
is not yet known.

Our starting point for constructing distorted subgroups of hyperbolic groups
is the following Lemma of Bestvina, Feighn and Handel \cite{BFH}:

\begin{lemma} \cite{BFH}
There exists a hyperbolic group $G$ such that $1 \rightarrow F \rightarrow G
\rightarrow F \rightarrow 1$ is exact, where $F$ is free of rank 3.
\label{BFH}
\end{lemma}

Let $F_1 \subset G$ denote the normal subgroup. Let $F_2 \subset G$ denote
a section of the quotient group. Let $G_1,{\cdots},G_n$ be $n$ distinct
copies of $G$. Let $F_{i1}$ and $F_{i2}$ denote copies of $F_1$ and $F_2$
respectively in $G_i$. Let \\ \begin{center}
$G = {G_1}{*_{H_1}}{G_2}*{\cdots}{*_{H_{n-1}}}{G_n}$
\end{center}
where each $H_i$ is a free group of rank 3, the image of $H_i$ in $G_i$
is $F_{i2}$ and the image of $H_i$ in $G_{i+1}$ is $F_{(i+1)1}$. Then $G$
is hyperbolic. This follows inductively from the fact that the image of 
$H_i$ in $G_i$ is quasiconvex in 
${G_1}{*_{H_1}}{G_2}*{\cdots}{*_{H_{i-1}}}{G_i}$ and the main combination theorem
of \cite{BF}.

Let $H = F_{11} \subset G$. Then the distortion of $H$ is superexponential
for $n > 1$. In fact, it can be readily checked that the distortion function 
is an iterated exponential of height $n$.

Note further that ${G_1}{*_{H_1}}{G_2}$ can be regarded as a graph of groups
with one vertex and three edges, where the vertex group is $G_1$ and edge groups
are isomorphic to $F$. Then from Corollary  \ref{main1}, the pair 
$({{G_1},{G_1}{*_{H_1}}{G_2}})$ 
has a Cannon-Thurston map. Proceeding inductively
and observing that a composition of Cannon-Thurston maps is a Cannon-Thurston
map, we see that $(H,G)$ has a Cannon-Thurston map. 

The next class of examples are not known to have Cannon-Thurston maps:

Our starting point is again Lemma \ref{BFH}. Let ${a_1},{a_2},a_3$ be 
generators of $F_1$ and ${b_1},{b_2},{b_3}$ be generators of $F_2$. Then
\begin{center}
$G = \{ {a_1},{a_2},{a_3},{b_1},{b_2},{b_3} : {b_i^{-1}}{a_j}{b_i} = w_{ij} \}$
\end{center}
where $w_{ij}$ are words in $a_i$'s. We add a letter $c$ conjugating $a_i$'s
to `sufficiently random' words in $b_j$'s to get $G_1$. Thus,
\begin{center}
$G_1 = \{ {a_1},{a_2},{a_3},{b_1},{b_2},{b_3},c : {b_i^{-1}}{a_j}{b_i} = w_{ij}, {c^{-1}}{a_i}c = v_i \}$,
\end{center}
where $v_i$'s are words in $b_j$'s satisfying a small-cancellation type condition to ensure that $G_1$ is hyperbolic. See \cite{Gromov}, pg. 151 for details
on addition of `random' relations.

It can be checked that these examples have distortion function greater than
any iterated exponential.

The above set of examples were motivated largely by examples of distorted
cyclic subgroups in \cite{Gromov2}, pg. 67.

So far, there is no satisfactory way of manufacturing examples of hyperbolic
subgroups of hyperbolic groups exhibiting arbitrarily high distortion.
It is easy to see that a subgroup of sub-exponential distortion is
quasiconvex \cite{Gromov2}. Not much else is known. For instance, one
does not know if $A^{n^2}$ can appear as a distortion function. 

The situation is far more satisfactory in the case of Kleinian groups. The 
following class of examples appears in work of Minsky \cite{Minsky1}:

Let $S$ be a hyperbolic punctured torus so that the two shortest geodesics
$a$ and $b$ are orthogonal and of equal length. Let $S_0$ denote $S$ minus
a neighborhood of the cusp. Let $N_{\delta}(a)$ and
$N_{\delta}(b)$ be regular collar neighborhoods of $a$ and $b$ in $S_0$.
For $n\in{\Bbb{N}}$, define ${\gamma_n} = a$ if $n$ is even and equal
to $b$ if $n$ is odd. Let $T_n$ be the open solid torus neighborhood
of ${\gamma_n}{\times}\{n+{\frac{1}{2}}\}$ in ${S_0}{\times}[0,{\infty})$
given by 
\begin{center}
$T_n = N_{\delta}({\gamma_n}){\times}(n,n+1)$
\end{center}
and let
$M_0 = ({S_0}){\times}[0,{\infty}){\setminus}{\bigcup}_{n\in{\Bbb{N}}}T_n$.

Let $a(n)$ be a sequence of positive integers greater than one. Let
${\hat{\gamma_n}} = {\gamma_n}{\times}\{{n}\}$ and let $\mu_n$ be an
oriented  meridian for $\partial{T_n}$ with a single positive intersection
with ${\hat{\gamma_n}}$. Let $M$ denote the result of gluing to each 
$\partial{T_n}$  a solid torus $\hat{T_n}$,
such that the curve ${{\hat{\gamma_n}}^{a(n)}}{\mu_n}$ is glued to a 
meridian.  Let
$q_{nm}$ be the mapping class 
from $S_0$ to itself obtained by identifying $S_0$ to ${S_0}{\times}m$,
pushing through $M$ to ${S_0}{\times}n$ and back to $S_0$. Then
$q_{n(n+1)}$ is given by $\Phi_n = 
D_{\gamma_n}^{a(n)}$, where $D_c^k$ denotes Dehn 
twist along $c$, $k$ times. Matrix representations of $\Phi_n$ are
given by \\ 
\[ \Phi_{2n} = \left( \begin{array}{cc}
                         1 & a(2n) \\ 0 & 1
                   \end{array} \right)  \]
and \[ \Phi_{2n+1} = \left( \begin{array}{cc}
                         1 & 0 \\ a(2n+1) & 1
                   \end{array} \right)  \]

Recall that the metric on $M_0$ is the restriction of the product metric.
$\hat{T_n}$'s are given hyperbolic metrics such that their boundaries
are uniformly quasi-isometric to $\partial{T_n} \subset M_0$. Then from
\cite{Minsky1}, $M$ is quasi-isometric to the complement of a rank
one cusp in the convex core of 
a hyperbolic manifold $M_1 = {{\Bbb{H}}^3}/{\Gamma}$.
 Let $\sigma_n$ denote the shortest path
from $S_0{\times}1$ 
to $S_0{\times}n$. Let ${\overline{\sigma_n}}$ denote $\sigma_n$ with reversed
orientation. Then $\tau_n = \sigma_n{\gamma_n}{\overline{\sigma_n}}$
is a closed path in $M$ of length $2n+1$. 
Further $\tau_n$ is homotopic to a curve 
$\rho_n = \Phi_1{\cdots}\Phi_n{(\gamma_n)}$ on $S_0$. Then
\begin{center}
$\Pi_{i=1{\cdots}n}{a(i)} \leq l({\rho_n}) \leq \Pi_{i=1{\cdots}n}{(a(i)+2)}$
\end{center}

Hence
\begin{center}
$\Pi_{i=1{\cdots}n}{a(i)} \leq (2n+1)disto(2n+1) \leq
  \Pi_{i=1{\cdots}n}{(a(i)+2)}$
\end{center}

Since $M$ is quasi-isometric to the complement of the cusp of a hyperbolic
manifold and $\gamma_n$'s lie in a complement of the cusp, the distortion
function of $\Gamma$ is of the same order as the distortion function above.
In particular, functions of arbitraily fast growth may be realised. This 
answers a question posed by Gromov \cite{Gromov2} pg.   66.

Manifolds with unbounded $a(n)$'s are not  known to have Cannon-Thurston
maps.

\medskip

{\bf Acknowledgements:} The author would like to thank his advisor Andrew
Casson for helpful comments and Curt Mcmullen for pointing out the results
of \cite{and-mask}.

\bibliography{tree}
\bibliographystyle{plain}

Address : Department of Mathematics, University of California, Berkeley, CA 94720, USA.

email : mitra@@math.berkeley.edu

\end{document}